\documentclass[12pt]{article}

\title{Geodesics in stationary spacetimes. Application to Kerr spacetime.}

\author{Jos\'{e} Luis Flores and Miguel S\'anchez \thanks{Research partially
supported by  
MCyT grant number BFM2001-2871-C04-01.} \\ Depto. Geometr\'{\i}a y
Topolog\'{\i}a, Fac. Ciencias, Univ. Granada, \\ Avda. Fuentenueva s/n, 18071-Granada, Spain.}

\usepackage{graphicx}
\usepackage{amssymb}
\usepackage{latexsym}
\usepackage[latin1]{inputenc}

\newfont{\bb}{msbm10 at 12pt}
\newfont{\bt}{msbm10 at 8pt}
\newfont{\btt}{msbm10 at 5pt}
\newfont{\kt}{cmti10 at 6pt}
\newfont{\bg}{msbm10 at 16pt}

\newtheorem{lemma}{Lemma}
\newtheorem{remark}{Remark}
\newtheorem{theorem}{Theorem}

\newtheorem{corollary}{Corollary}

\font\ddpp=msbm10 scaled \magstep 1 %Caracteres "doble palo".
\def\R{\hbox{\ddpp R}}    % numeros reales
\def\N{\hbox{\ddpp N}}    % números naturales
\def\D{\hbox{\ddpp D}}
\def\P{\hbox{\ddpp P}}              

\def\be{\begin{equation}}
\def\ee{\end{equation}}

\begin{document}
\maketitle

%\begin{center}
%{\it Running title:} GRW spacetimes
%\end{center}

\begin{abstract}

The Levi-Civita connection and geodesic equations for a stationary spacetime are studied in depth. General formulae which generalize those for warped products are obtained. These results are applicated to some regions of Kerr spacetime previously studied by using variational methods. We show that they are neither space-convex nor geodesically connected. Moreover, the whole stationary part of Kerr spacetime is not geodesically connected, except when the angular momentum is equal to zero (Schwarzschild spacetime).

\end{abstract}

\newpage 
\section{Introduction}

The existence of a Killing vector field $K$ on a spacetime $(M,<\cdot, \cdot >)$ is specially useful to study its geometry. It is well--known that, around each point $p$ such that 
$K_p \neq 0$, coordinates
$(t, x^1, \dots , x^n)$ can be chosen such that $K = \partial_{t}$ and all the components 
$g_{ij}$ of the metric are independent of $t$; this justifies the name ``stationary'' for the spacetime when $K $ is timelike. The stationary observers along $K$ not only see a non-changing metric but also find a constant  $E= <\partial_t , \gamma' >$ for any geodesic $\gamma$; thus, photons and freely falling particles has constant {\em energy} $E$ for these observers.
When a non-vanishing Killing vector field $K$ is irrotational, i.e. the orthogonal distribution $K^{\perp}$ is involutive, then a local warped product structure appears ($g_{ti}= 0$, for $i=1, \dots, n$); if, additionally, $K$ is  timelike ({\em static} case), the observers along $K$
measure a metric with no cross terms between space and time. 

Since the introduction by Bishop and O'Neill of Riemannian warped products \cite{BO}, warped 
structures has been widely studied. This includes the Lorentzian case, where 
the contributions by O'Neill 
\cite{O-83} and Beem, Ehrlich and Powell \cite{BEP} (see also \cite{BEE})  have been especially relevant. 
Recall that there are many examples of warped products among classical relativistic spacetimes. Nevertheless, the non-warped stationary case becomes  difficult (we refer to \cite{Sa} for a summary of mathematical properties in this general case).

Among the more interesting problems in stationary spacetimes appear those related to geodesics, as geodesic completeness or geodesic connectedness. The former has been widely studied by Romero and one of the authors \cite{RSa}, \cite{RSb}. The latter was studied from a variational viewpoint first by Benci, Fortunato and Giannoni \cite{BF}, \cite{BFG}, and then by several authors 
(see the book \cite{Ma}). The more classical stationary spacetime is  Kerr spacetime which, essentially, represents the stationary  gravitational field outside a rotating star. This spacetime is characterized by the mass $m>0$ and angular momentum $ma$ of the star; in the limit case $a=0$ it becomes Schwarzschild (outer) spacetime, which is static. 
The more recent book by O'Neill \cite{O} study especifically Kerr spacetime (including the extending  ---non-stationary--- regions). 

The aim of this article is to study geodesics in a stationary standard spacetime
$M=\R \times M_0$, giving some applications to Kerr spacetime. It is organised as follows. 

In Section 2 we give general formulas for the Levi-Civita connection, which
 are an extension of those for warped products. In fact, they remain valid 
 if $K=\partial_t$ is not timelike, and they 
can be generalized to more general warped products with crossed terms (Theorem \ref{t1}, 
Remark \ref{rem1}). Then, we study geodesic equations, 
generalizing those in \cite{Sa-non} (Theorem \ref{t2}, Remark \ref{rem3}).
Finally, we give general formulas for the Hessian of a function 
$\phi$ independent of $t$, Hess$\phi$ (Theorem \ref{t3}). This Hessian 
is directly involved in the problem of geodesic connectedness because of
the following notion of convexity. Let $c$ be a regular value of $\phi$ and 
$D=\phi^{-1}((c,\infty))$. The boundary $\partial D = \phi^{-1}(c)$ is 
(time, light, space-like) {\it convex} if Hess$\phi (\tilde{v},\tilde{v})\leq 0$ for all
(time, light, space-like) vectors $\tilde{v}$ tangent to $\partial D$. Under certain
natural assumptions, the convexity of the boundary imply geodesic 
connectedness  (see \cite{Ma} or \cite{Sa-non2}).

In Section 3 we study some general properties of 
(the stationary part of) Kerr spacetime $M^a$.
The slow rotating ($a^2 < m^2$) or limit ($a^{2}=m^{2}$) cases are specially interesting because of the properties of the boundary of the stationary part; 
moreover, Schwarzschild spacetime is also included 
as the limiting static case $a=0$. Section \ref{s3} is centered on the cases $a^{2}\leq m^{2}$; however, our study will also cover the fast rotating case $a^{2}>m^{2}$ (Theorem \ref{te}). In particular, we study regions
$M^{a}_{\epsilon}$ introduced by Giannoni and Masiello \cite{GM},
\cite{Ma}. Recall that, from a variational viewpoint, 
the boundary of $M^a$ is singular and difficult to study,
even in the static case \cite{bfg2}. Thus, regions $M^{a}_{\epsilon}$,
with smooth boundary but arbitrarily close to the boundary of
$M^{a}$, were introduced by these authors. For $a$ and $\epsilon$
small enough, they showed that the boundary of $M^{a}_{\epsilon}$
is time and light convex, which allowed them to prove some results 
on the existence of causal geodesics. In Section 3 general 
formulas to study convexity are provided, and $M^{a}_{\epsilon}$
is shown to be non-space convex, Corollary \ref{cec}. 

In Section 4 we use geodesic equations to prove directly that neither
Kerr spacetime $M^a$ for any $a\neq 0$ nor any region $M^{a}_{\epsilon}$, $\epsilon >0$
(including the case $a=0$) are geodesically connected.

Finally, in Section 5 we prove that the excluded case 
($a=0=\epsilon$, Schwarzschild spacetime) is geodesically connected. 
Even though the geodesic connectedness of this spacetime 
has already been proven by using variational methods 
in \cite{bfg2}, we include this proof because of several reasons:
(1) it is completely different, based on topological arguments 
introduced by the authors in \cite{FS}, and fully
adapted to this case, (2) it is easily translatable to 
 Schwarzschild black hole (Remark \ref{rr}), where variational methods  
seems to fail, and (3) it is a simple case of the more  
involved proof of the geodesic connectedness of outer Kerr spacetimes
(where the causal character of $\partial_{t}$ changes, if 
$a\neq 0$), which will be the subject of a next article \cite{FS-pr}.

\section{Levi-Civita connection and geodesics} \label{s2}

Let $ (M_{0},<.,.>_{R}) $ be a Riemannian manifold, $ (\R,-dt^{2})$  the set of the real numbers with its usual metric reversed, and $\delta$ and $\beta$ a vector field and a positive smooth function on $M_{0}$, respectively. A {\it (standard) stationary spacetime} is the product manifold $M=\R \times M_{0}$ endowed with the Lorentz metric
\begin{equation}
\label{prim}
<.,.>=-\beta(x)dt^{2}+<.,.>_{R}+2<\delta,.>_{R}dt.
\end{equation}

Choosing an orthonormal basis $B_{0}=(e_{1},\ldots,e_{n})$ at some $T_{x}M_{0}$ the matrix of $<.,.>$ in $B=(\partial_{t},e_{1},\ldots, e_{n})$ is
\begin{equation}\label{e-10}
\left( \begin{array}{c|ccc} -\beta & \delta^{1} & \ldots & \delta^{n} \\ \hline \delta^{1} & 1 &  & 0 \\ \vdots & & \ddots & \\ \delta^{n} & 0 &  & 1 \end{array} \right) \equiv        \left( \begin{array}{c|c} -\beta & \delta^{t} \\ \hline \delta & Id_{n} \end{array} \right)
\end{equation}
where $\delta^t = (\delta^{1}, \cdots , \delta^{n})$  are the components of $\delta$ in $B_{0}$, and $Id_{n}$ is the identity matrix $n\times n$. Putting $\Lambda =-1/(\beta +\parallel \delta \parallel_{R}^{2})$, the inverse of (\ref{e-10}) is
\begin{equation}\label{e0}
 \left( \begin{array}{c|c} \Lambda & -\Lambda \delta^{t} \\ \hline  &  \\ -\Lambda\delta & Id_{n}+\Lambda\delta\otimes\delta^{t} \\ & \end{array} \right)
\end{equation}
where $\parallel .\parallel_{R}$ denotes the $<.,.>_{R}-\mbox{norm}$.

From now on, $\nabla$ and $\nabla^{R}$ will denote the Levi-Civita connection of $<.,.>$ and $<.,.>_{R}$, respectively. For each vector field $V$ on $M_{0}$, $V\in \mathfrak{X} (M_{0})$, its lifting to $M$ will be denoted $\overline{V}$: that is, $\overline{V}_{(t,x)}=V_{x}$, $\forall (t,x)\in M$ (analogously, if necessary, for a vector field on $\R$).

\begin{theorem}\label{t1} Let $(\R\times M_{0},<.,.>)$ be a stationary spacetime and let $V,W\in \mathfrak{X} (M_{0})$. Then:

(i)
\begin{equation}\label{e2}
\nabla_{\partial_{t}}\partial_{t}=-\frac{1}{2}\Lambda<\delta,\nabla^{R}\beta>_{R}\partial_{t}+\frac{1}{2}\Lambda<\delta,\nabla^{R}\beta>_{R}\delta +\frac{1}{2}\nabla^{R}\beta.
\end{equation}

(ii)
\begin{equation}\label{e1}
\begin{array}{ll}
\nabla_{\overline{V}}{\overline{W}}= & \frac{1}{2}\Lambda(<W,\nabla_{V}^{R}\delta>_{R}+<V,\nabla_{W}^{R}\delta>_{R})\partial_{t} \\ & -\frac{1}{2}\Lambda(<W,\nabla_{V}^{R}\delta>_{R}+<V,\nabla_{W}^{R}\delta>_{R})\delta+\nabla_{V}^{R}W.
\end{array}
\end{equation}

(iii)
\begin{equation}\label{e3}
\begin{array}{ll}
2\nabla_{\overline{V}}\partial_{t}=2\nabla_{\partial_{t}}\overline{V}= & -\Lambda(V(\beta)+<\delta,\nabla_{V}^{R}\delta>_{R}-<\nabla_{\delta}^{R}\delta,V>_{R})\partial_{t} \\ & +\Lambda(V(\beta)+<\delta,\nabla_{V}^{R}\delta>_{R}-<\nabla_{\delta}^{R}\delta,V>_{R})\delta \\ & +\nabla_{V}^{R}\delta -<\nabla^{R}_{(\cdot)}\delta,V>_{R}^{\natural}
\end{array}
\end{equation}
where $^{\natural}$ denotes the vector field on $M_{0}$ 
metrically associated to the corresponding 1-form 
(that is, 
 $<Y,<\nabla^{R}_{(\cdot)}\delta,V>_{R}^{\natural}>_{R}= <\nabla_{Y}^{R}\delta,V>_{R}$ for any $Y\in \mathfrak{X}(M_{0})$).

\end{theorem}

{\it Proof.} First, recall Koszul's formula
\begin{equation}\label{e-100}
\begin{array}{ll}
2<\nabla_{Y}Z,X>= & Y<Z,X>+Z<X,Y>-X<Y,Z> \\ & -<Y,[Z,X]>+<Z,[X,Y]>+<X,[Y,Z]>
\end{array}
\end{equation}
for any vector fields $Y,Z,X \in \mathfrak{X}(M)$.

(i) Clearly, $2<\partial_{t},\nabla_{\partial_{t}}\partial_{t}>=\partial_{t}
<\partial_{t},\partial_{t}>=0$. For any $X\in\mathfrak{X} (M_{0})$, as $\partial_{t}<\partial_{t},\overline{X}>=\partial_{t}<\delta,X>_{R}=0$, thus:

\[
\begin{array}{ll}
2<\nabla_{\partial_{t}}\partial_{t},\overline{X}>= & 2\partial_{t}<\partial_{t},\overline{X}>-\overline{X}<\partial_{t},\partial_{t}>= \\ & =X(\beta)=<\nabla^{R}\beta,X>_{R}.
\end{array}
\]
Therefore, the result follows multiplying the inverse matrix (\ref{e0}) by the components 
$(0, \frac{1}{2}\nabla^{R}\beta )$ of the 1-form associated to $\nabla_{\partial_{t}}\partial_{t}$.

%\[
%\left( \begin{array}{cc} a & -a\delta^{t} \\  &  \\ -a\delta & Id_{n}+a\delta\otimes\delta^{t} \\ & \end{array} \right)
%\cdot \left( \begin{array}{c} 0 \\ \frac{1}{2}\nabla^{R}\beta \end{array} \right)=\left( \begin{array}{c} -a\frac{1}{2}<\delta,\nabla^{R}\beta> \\ \frac{1}{2}\nabla^{R}\beta +a\frac{1}{2}<\delta,\nabla^{R}\beta>\delta \end{array}\right)
%\]

%that is 

 %Finally, consider the vector field $\nabla_{\partial_{t}}V=\nabla_{V}\partial_{t}$ with $V$ tangent to $M_{0}$, then

(ii) Fixed $x_{0}\in M$ there is no loss of generality if we assume $[V,W]_{x_{0}}=0$. Let               $X\in \mathfrak{X} (M_{0})$ satisfying  $[X,W]=[X,V]=0$ at $x_{0}$. So, we have at this point by using (\ref{e-100}): 
\[
\begin{array}{ll}
2<\nabla_{\overline{V}}\overline{W},\overline{X}>= & V<W,X>_{R}+W<X,V>_{R}-X<V,W>_{R}= \\ & =2<\nabla_{V}^{R}W,X>_{R}
\end{array}
\]
Moreover, using $\nabla_{V}^{R}W+\nabla_{W}^{R}V=2\nabla_{V}^{R}W$ at $x_{0}$ and again by (\ref{e-100}),

\[
\begin{array}{ll}
2<\nabla_{\overline{V}}\overline{W},\partial_{t}>= & V<W,\delta >_{R}+W<V,\delta >_{R}= \\ & =<W,\nabla_{V}^{R}\delta >_{R}+<V,\nabla_{W}^{R}\delta >_{R}+2<\delta,\nabla_{V}^{R}W>_{R}.
\end{array}
\]
Finally, the components of $\nabla_{\overline{V}}\overline{W}$ are obtained again by using (\ref{e0}).
%\[
%\left( \begin{array}{cc} a & -a\delta^{t} \\  &  \\ -a\delta & Id_{n}+a\delta\otimes\delta^{t} \\ & \end{array} \right)
%\left( \begin{array}{c} <\delta,\nabla_{V}^{R}W+\frac{1}{2}\{ <W,\nabla_{V}^{R}\delta>_{R}+<V,\nabla_{W}^{R}\delta>_{R}\} \\ \nabla_{V}^{R}W \end{array} \right)
%\]

%that is

 %Consider now the vector field $\partial_{t}$, then we have for  $\nabla_{\partial_{t}}\partial_{t}$

(iii) Clearly 
\[
2<\nabla_{\overline{V}}\partial_{t},\partial_{t}>=-V(\beta)=-<\nabla^{R}\beta,V>_{R}.
\]
Taking $X\in \mathfrak{X}(M_{0})$ commuting with $V$ at $x_{0}$, we have from (\ref{e-100})

\[
\begin{array}{ll}
2<\nabla_{\overline{V}}\partial_{t},\overline{X}>= & \overline{V}<\partial_{t},\overline{X}>+\partial_{t}<\overline{V},\overline{X}>-\overline{X}<\overline{V},\partial_{t}>= \\ & =V<\delta,X>_{R}-X<V,\delta >_{R}=<\nabla_{V}^{R}\delta,X>_{R} \\ & -<\nabla_{X}^{R}\delta,V>_{R},
\end{array}
\]
and the result follows by using (\ref{e0}) again. $\Box$
%\[ 
%\begin{array}{l}  \left( \begin{array}{cc} a & -a\delta^{t} \\  &  \\ -a\delta & Id_{n}+a\delta\otimes\delta^{t} \\ & \end{array} \right)\cdot \frac{1}{2} \left( \begin{array}{c} -V(\beta) \\ \nabla_{V}^{R}\delta-<\nabla^{R}\delta,V>_{R}^{\natural} \end{array} \right)= \\ =\frac{1}{2} \left( \begin{array}{c} -aV(\beta)-a<\delta,\nabla_{V}^{R}\delta>_{R}+a<\nabla_{\delta}^{R}\delta,V>_{R} \\ a\{V(\beta)+<\delta,\nabla_{V}^{R}\delta>_{R}-<\nabla_{\delta}^{R}\delta,V>_{R}\}\delta +\nabla_{V}^{R}\delta -<\nabla^{R}\delta,V>_{R}^{\natural} \end{array} \right) \end{array}
%\] 

\begin{remark}\label{rem1}
\em{Let $\nabla^{R}\delta(V,W) = <\nabla_{V}^{R}\delta,W>_{R}$, and consider its symmetric and skew-symmetric parts: Sym$\nabla^{R}\delta(V,W) = (\nabla^{R}\delta(V,W) + \nabla^{R}\delta(W,V))/2$; 
Sk$\nabla^{R}\delta(V,W) = (\nabla^{R}\delta(V,W) - \nabla^{R}\delta(W,V))/2$. In what follows,  rot$\delta$ = 2Sk$\nabla^{R}\delta $.  In Theorem \ref{t1}, the differences between the general stationary case and the static (warped) case are the following:

(i) In (\ref{e2}), the term:  
$$-\frac{1}{2}\Lambda<\delta,\nabla^{R}\beta>_{R}
(\partial_{t} - \delta ).$$ 
Note that this term vanishes if $\beta$ is constant. This happens for 
the conformal metric $<\cdot , \cdot >/\beta$, which can be used to study null geodesics
(recall that null geodesics are conformal-invariant, up to reparametrizations).

(ii) In (\ref{e1}), the term: 
$$\Lambda {\rm Sym}\nabla ^{R}\delta (V,W) (\partial_t - \delta).$$
This term vanishes if $\nabla^R \delta$ is skew symmetric, that is, if $\delta $ is a Killing vector field. 

(iii) In (\ref{e3}), the terms:  $$  
\Lambda V(\beta) \delta-\Lambda {\rm rot}\delta (V,\delta) (\partial_t - \delta) 
+{\rm rot}\delta (V, \cdot )^{\natural}.$$
These terms reduces to $\Lambda V(\beta) \delta$
 if  $\delta $ is irrotational, that is, locally a gradient vector field. 
}
\end{remark}

For geodesic equations, let  $\gamma(s)=(t(s),x(s))$ be a geodesic on $M$ with initial condition $\gamma'(0)=(t'(0),x'(0))$ and $t'(0)\neq 0\neq x'(0)$ (the modification otherwise would be straightforward). We can choose a vector field $X$ (resp. $T$) on $ M_{0}$ (resp. $\R$) extending $x'(s)$ (resp. $t'(s)$). Then, the vector field $Z=\overline{T}+\overline{X}$ satisfies on $\gamma$:
% geodesics with initial condition $(t'(0),\hat{X}(x))$ on the hypersurface $t=t(0)$ yield a vector field $Z=(T,X)$ around $\gamma(0)$ satisfying $\nabla_{Z}Z=0$, being the geodesic $\gamma$ an integral curve of $Z$. So, we have
\begin{equation}\label{e4}0=\nabla_{Z}Z=\nabla_{\overline{T}}\overline{T}+\nabla_{\overline{T}}\overline{X}+\nabla_{\overline{X}}\overline{T}+\nabla_{\overline{X}}\overline{X} \end{equation}
and we can use the equalities (\ref{e2}), (\ref{e1})  and (\ref{e3}) to rewrite this relation. 

Note that on the geodesic $\gamma$

\begin{equation}\label{e5} \left\{ \begin{array}{l} X(x(s))=x'(s) \\ T(t(s))=t'(s)\partial_{t} \end{array} \right. \end{equation} 
holds, and so

\begin{equation}\label{e6} \left\{ \begin{array}{l} \nabla_{\overline{T}}\overline{T}=t'\partial_{t}(t')\partial_{t}+t'^{2}\nabla_{\partial_{t}}\partial_{t} \\ \nabla_{\overline{T}}\overline{X}=\nabla_{\overline{X}}\overline{T}=t'\nabla_{\partial_{t}}\overline{X}. \end{array} \right. \end{equation}

In order to obtain an equation for $x(s)$ we will use that, in the last member of (\ref{e4}), the sum of the components on $M_{0}$ of the four vector fields must be $0$. These components can be obtained from Theorem \ref{t1}. Using also (\ref{e6}) and writing $W^{M_{0}}\equiv$ the projection of the vector field $W$ on $M_{0}$,
\[
\begin{array}{ll}(\nabla_{\overline{T}}\overline{T})^{M_{0}}= & t'^{2}(\nabla_{\partial_{t}}\partial_{t})^{M_{0}}= \\ & =\frac{1}{2}t'^{2}(\Lambda<\delta,\nabla^{R}\beta>_{R}\delta+\nabla^{R}\beta) \\ (\nabla_{\overline{T}}\overline{X}+\nabla_{\overline{X}}\overline{T})^{M_{0}}= & 2t'(\nabla_{\overline{X}}\partial_{t})^{M_{0}}= \\ & =\Lambda t'\left(X(\beta)
+{\rm rot}\delta(X,\delta)\right)\delta 
%\\ & 
+t' {\rm rot}\delta(X,\cdot)^{\natural}  \\
(\nabla_{\overline{X}}\overline{X})^{M_{0}}= & \nabla_{X}^{R}X-\Lambda \nabla^{R}\delta(X,X)\delta.
\end{array}
\]

Adding these three relations and composing with $\gamma$, we have

\begin{equation}\label{e7}
\begin{array}{ll} 
\nabla_{x'}^{R}x'= & \Lambda \nabla^{R}\delta(x',x') \delta 
\\ & -\Lambda t' \left(<\nabla^{R}\beta,x'>_{R} +{\rm rot}\delta(x',\delta)\right)\delta 
\\ & -t' {\rm rot}\delta(x',\cdot )^{\natural} 
\\ & -\frac{1}{2}t'^{2}(\Lambda <\delta,\nabla^{R}\beta>_{R}\delta +\nabla^{R}\beta )
\end{array}
\end{equation}
On the other hand, %recall that%
  $q=<\gamma ',\gamma '> $  is constant for any geodesic $\gamma$ and,  as $\partial_{t}$ is a Killing vector field,  $E =<\gamma ',\partial_{t}>$ is constant too. That is, we also have the relations: 

\begin{equation}\label{e8} \left\{ \begin{array}{l} <(t',x'),(t',x')>=-\beta t'^{2}+2<\delta,x'>_{R}t'+<x',x'>_{R}=q \\ <\partial_{t},(t',x')>=-\beta t'+<\delta,x'>_{R}=E.
\end{array} \right. \end{equation} 
Let $\mathfrak{X}_{r,s}(M_0)$ be the space of $r-$contravariant, $s-$covariant tensor fields on $M_0$ ($\mathfrak{X} (M_0)\equiv \mathfrak{X}_{1,0} (M_0)$; 
$\mathfrak{X}_{1,1}(M_{0})$ is identifiable to the space of endomorphism fields).
Equation (\ref{e7}) can be written as 

\begin{equation}\label{e9} \nabla_{x'}^{R}x'=t'^{2}R_{0}(x)+t'R_{1}(x,x')+R_{2}(x,(x',x'))
\end{equation}
with $R_0 \in \mathfrak{X}(M_{0}), R_1 \in \mathfrak{X}_{1,1}(M_{0}), R_2 \in 
\mathfrak{X}_{1,2}(M_{0})$
putting:
\begin{equation} \label{eee} \left\{
\begin{array}{l}R_{0}(x)=-\frac{1}{2}\left( \Lambda <\delta,\nabla^{R}\beta>_{R}\delta+\nabla^{R}\beta \right)(x) \\
R_{1}(x,x')= -\Lambda \left(<\nabla^{R}\beta,x'>_{R}+ 
{\rm rot}\delta (x',\delta)\right) \delta(x)+
{\rm rot}\delta (\cdot,x')^{\natural}  \\
R_{2}(x,(x',\overline{x}'))= \Lambda {\rm Sym}\nabla^{R}\delta(x',\overline{x}')\delta(x).
\end{array} \right.
\end{equation}

%\[
%\begin{array}{l}
%\psi(x')=-<\nabla^{R}\beta,x'>_{R}+<\nabla_{\delta}^{R},x'>_{R}-
%<\delta,\nabla_{x'}^{R}\delta>_{R}\in \mathfrak{X}_{0,1}(M_{0}) \\
% F_{\delta}(x')=<\nabla^{R}\delta,x'>_{R}^{\natural}-\nabla_{x'}^{R}\delta \in \mathfrak{X}_%%%{1,1}(M_{0})\quad (skew-adjoint)
% \\
%\nabla^{R}\delta(x',x')=<\nabla_{x'}^{R}\delta,x'>\in \mathfrak{X}_{2,0}(M_{0}).
%\end{array}
%\]
\begin{remark}
{\em 
$R_0$ vanishes if $\beta$ is constant, $R_1$ reduces to 
$-\Lambda <\nabla^{R}\beta,x'>_{R} \delta(x)$
 if $\delta$ is irrotational, and $R_2$ vanishes if $\delta$ is Killing.
}
\end{remark}
Substituting in (\ref{e9}) the value of $t'$ from the second equation in (\ref{e8}) a second order equation  for the spacelike component $x(s)$  is obtained. Then, the first relation (\ref{e8}) can be regarded as a first integral of this equation. 
 Moreover, a geodesic can be reconstructed for any solution of this differential equation.
Summing up:

\begin{theorem}\label{t2} Consider a curve $\gamma (s)=(t(s),x(s))$ in a stationary spacetime $(\R\times M_{0},<.,.>)$. The curve $\gamma$ is a geodesic if and only if 
 $E=<\gamma',\partial_{t}>$ is a constant and, then, $x(s)$ satisfy
\begin{equation}\label{e10} \nabla_{x'}^{R}x'=\overline{R}_{0}(x)+\overline{R}_{1}(x,x')+\overline{R}_{2}(x,x'\otimes x')
\end{equation}
where: %putting $E_{\beta}\equiv \frac{E}{\beta}$: $\delta_{\beta}\equiv\frac{\delta}{\beta}$
\[
\begin{array}{ll}
\overline{R}_{0}(x)= &  \frac{E^2}{\beta^2} R_{0}(x)  \\
 %& (=-\frac{E^{2}}{2\beta^{2}}\{\Lambda <\delta,\nabla^{R}\beta>\delta+\nabla^{R}\beta\}=\frac{\lambda^{2}}{2}\{\Lambda<\delta,\nabla^{R}\frac{1}{\beta}>\delta+\nabla^{R}\frac{1}{\beta}\})
\overline{R}_{1}(x,x')= & -\frac{E}{\beta}\left( 2\frac{<\delta,x'>_R}{\beta}R_{0}(x)+
R_{1}(x,x')\right) \\ 
%& (=\frac{E}{\beta}\{<\frac{\delta}{\beta},x'><\delta,\nabla^{R}\beta>\Lambda \delta+<\frac{\delta}{\beta},x'>\nabla^{R}\beta\}-\frac{E}{\beta}\{\psi(x')\delta(x)+F_{\delta}(x')\}= \\ & 
% =\frac{E^{2}}{2}\{\Lambda <\delta,\nabla^{R}\frac{1}{\beta}>\delta+\nabla^{R}\frac{1}{\beta}\})\\ 
\overline{R}_{2}(x,x'\otimes x')= & \frac{(<\delta,x'>_{R})^2}{\beta^{2}}R_{0}(x)+
\frac{<\delta,x'>_R}{\beta}R_{1}(x,x')+R_{2}(x,(x',x')),
% \\ & 
%(=\frac{1}{2\beta}<\delta,x'>^{2}\{\Lambda <\delta,\nabla^{R}\frac{1}{\beta}>\delta+\nabla^{R}\frac{1}{\beta}\} \\ &  +\frac{1}{\beta}<\delta,x'>\{\Lambda \psi(x')\delta+F_{\delta}(x')\}+\Lambda \nabla^{R}\delta(x',x')\delta= \\ & =(\frac{1}{2\beta}<\delta,x'>^{2}<\delta,\nabla^{R}\frac{1}{\beta}>+\frac{1}{\beta}<\delta,x'>\psi(x')+\nabla^{R}\delta(x,x'))\Lambda \delta \\ &  +\frac{1}{2\beta}<\delta,x'>^{2}\nabla^{R}\frac{1}{\beta}+\frac{1}{\beta}<\delta,x'>F_{\delta}(x')).
\end{array}
\]
and $R_0, R_1, R_2$ are as in (\ref{eee}). 
\end{theorem}

\begin{remark} \label{rem3}
{\em This result extends the construction in \cite{Sa-non} for the static case. Indeed, for any solution $x(s)$ of (\ref{e10}), the second equation (\ref{e8}) yields the value of $t(s)$. Moreover, all the geodesics can be reparametrized in such a way that $E = 0,1$. Thus, we have the following two cases: 

(a) Case $E=0$. The geodesic $\gamma$ is always spacelike and orthogonal to $\partial_t$, and we can also assume that the value of $q$ is fixed equal to 1.
Recall that in this case $\overline{R}_{0}=\overline{R}_{1}=0$. In the static case, 
%(or, with more generality, if $\delta$ is a Killing vector field)
equation (\ref{e10}) is just the equation of the geodesics of $(M_0, <\cdot,\cdot>_R)$; thus, these geodesics can be regarded as trivial. Nevertheless, when $\delta $ is not null the differential equation becomes:
$$\nabla_{x'}^{R}x'=\overline{R}_{2}(x,x'\otimes x') ,$$
(compare, for example, with \cite{Ma-93}).

(b) Case $E=1$. Equation (\ref{e10}) becomes rather complicated in general. In the static case, 
$x(s)$ satisfies the equation of a classical Riemannian particle under the potential 
$V=-1/2\beta $, as studied in \cite{Sa-non}, and the constant $q/2$ is the classical energy (kinetic plus potential) of this particle.  In the general case this interpretation does not hold\footnote{This is a reason for the different notation of the constants with respect to 
\cite{Sa-non}: $(\lambda,\epsilon)$ in this reference becomes $(E,q)$ now.}, even though $q$ is a constant of the motion.}
\end{remark}

Now, let us study the expression of the Hessian for the stationary spacetime. Let $\phi$ be a smooth function defined on a stationary spacetime $M=\R\times M_{0}$ and let $\tilde{V}$ be a vector field on $M$. At any point $p=(t,x)\in M$  
\[
{\rm Hess}_{\phi}(p)[\tilde{V},\tilde{V}]=\tilde{V}\tilde{V}\phi(p)-\nabla_{\tilde{V}}\tilde{V}(\phi)(p).
\]
We can assume $\tilde{V}=\overline{T}+\overline{V}$ being $T$ and $V$ vector fields on $\R$ and $M_{0}$, respectively, and suppose $T_{p}=t'_0\partial_{t}$. Next,  if $\phi$ is independent of $t$; we obtain at $p$
\begin{equation}\label{e11}
\begin{array}{ll}
{\rm Hess}_{\phi}(p)[\tilde{V},\tilde{V}]= & \overline{V}\overline{V}\phi -\nabla_{\overline{V}}\overline{V}(\phi)-\nabla_{\overline{T}}\overline{V}(\phi)-\nabla_{\overline{V}}\overline{T}(\phi)-\nabla_{\overline{T}}\overline{T}(\phi)= \\  &  =\overline{V}\overline{V}\phi-\nabla_{\overline{V}}\overline{V}(\phi)-t'_0\nabla_{\partial_{t}}\overline{V}(\phi)-t'_0\nabla_{\overline{V}}\partial_{t}(\phi)-t'^{2}_{0}\nabla_{\partial_{t}}\partial_{t}(\phi).
\end{array}
\end{equation}  

Taking into account (\ref{e2}), (\ref{e1}) and (\ref{e3}) and using, again, that $\phi$ is independent of $t$ we have
\begin{equation}\label{e12}
\begin{array}{ll}
\nabla_{\overline{V}}\overline{W}(\phi)= & -\Lambda {\rm Sym}\nabla^R\delta(V,W) \delta (\phi)+\nabla_{V}^{R}W(\phi) \\ 
\nabla_{\partial_{t}}\partial_{t}(\phi)= & \frac{1}{2}\Lambda <\delta,\nabla^{R}\beta>_{R}\delta (\phi) +\frac{1}{2}\nabla^{R}\beta (\phi) 
\\ 2\nabla_{\overline{V}}\partial_{t}(\phi)= & \Lambda \left(V(\beta)+
{\rm rot}\delta(V,\delta)\right) \delta (\phi) %\\ & 
+{\rm rot}\delta (V, \nabla^R\phi) 
=2\nabla_{\partial_{t}}\overline{V}(\phi).
\end{array}
\end{equation}
Finally, if we substitute (\ref{e12}) in (\ref{e11}) we obtain:
\begin{theorem}\label{t3} 
Let $\phi : \R\times M_{0}\rightarrow \R$ be a function independent of $t$, and $\tilde{v}=t'_0\partial_{t}+v$ a vector on $T_{(t,x)}M=T_{t}\R\times T_{x}M_{0}$. Then
\begin{equation}\label{e13}
\begin{array}{ll}{\rm Hess}_{\phi}(p)[\tilde{v},\tilde{v}]= & {\rm Hess}_{\phi}^{R}(x)[v,v] \\ 
& -\Lambda \left( t'_{0}v(\beta) + \frac{1}{2} t'^{2}_{0}\delta(\beta) -{\rm Sym} \nabla^R \delta (v,v) +
t'_{0} {\rm rot}\delta(v,\delta) \right)\delta (\phi)
\\ & + t'_0{\rm rot}\delta(\nabla^{R}\phi,v) -\frac{1}{2}t'^{2}_{0}\nabla^{R}\beta(\phi).
\end{array}
\end{equation} 
\end{theorem}

\section{Kerr spacetime}\label{s3}

{\it Kerr spacetime} represents the stationary axis-symmetric asymptotically flat
gravitational field outside a rotating massive object. Let $m>0$ and $a$ be two constants, such that $m$ represents the mass of the object and $ma$ the angular momentum as measured from infinity. In the Boyer-Lindquist coordinates $(t,r,\theta,\varphi )$, Kerr metric takes the form
\begin{equation}\label{ep}
ds^{2}=g_{4,4}dt^{2}+g_{1,1}dr^{2}+g_{2,2}d\theta^{2}+g_{3,3}d\varphi^{2}+2g_{3,4}dt d\varphi 
\end{equation}
with
\begin{equation}
\begin{array}{lll}g_{1,1}=\frac{\lambda(r,\theta)}{\triangle (r)}\quad\quad & g_{3,3}=[r^{2}+a^{2}+\frac{2mra^{2}\sin^{2}\theta}{\lambda(r,\theta)}]\sin^{2}\theta\quad\quad & g_{4,4}=-1+\frac{2mr}{\lambda(r,\theta)} \\ g_{2,2}=\lambda(r,\theta)\quad\quad & g_{3,4}=-\frac{2mra\sin^{2}\theta}{\lambda(r,\theta)}\quad\quad &  
\end{array}
\end{equation}
and being
%\begin{equation}\label{ep}
%ds^{2}=\lambda (\frac{dr^{2}}{\triangle}+d\theta^{2})+(r^{2}+a^{2})\sin^{2}\theta d\varphi^{2}-dt^{2}+\frac{2mr}{\lambda}(a\sin^{2}\theta d\varphi -dt^{2})^{2},
%\end{equation}
\[
\lambda(r,\theta)=r^{2}+a^{2}\cos^{2}\theta\quad\hbox{and}\quad \triangle (r)=r^{2}-2mr+a^{2}.
\]
So, (\ref{ep}) can be written as in (\ref{prim}) taking
\begin{equation}
\begin{array}{l}
<.,.>_{R}=g_{1,1}dr^{2}+g_{2,2}d\theta^{2}+g_{3,3}d\varphi^{2} \\ \delta=\frac{g_{3,4}}{g_{3,3}}\partial_{\varphi} \\ \beta(r,\theta)=-g_{4,4}
\end{array}
\end{equation}
The structure of Kerr spacetime depends on the physical constants of the object $m$ and $a$.
In what follows, we will consider the case $a^{2} \leq m^{2}$ (the case $a^{2}>m^{2}$ is simpler and the conclusions for this case are summarized at the end of Section \ref{sad}). The function $\triangle (r)$ has the zeroes $r_{+}=m+\sqrt{m^{2}-a^{2}}$ and $r_{-}=m-\sqrt{m^{2}-a^{2}}$. The hypersurfaces $r=r_+ , r=r_-$
%$\R\times\{x\in \R^{3}:r=r_{+}\}\times S^^2$ and $\R\times\{x\in \R^{3}:r=r_{-}\}\times S^2$ 
are singular for (\ref{ep}); they are the {\it event horizons}. Outside the first one 
({\it outer Kerr spacetime}, $r>r_+$) the metric is not stationary if $a^{2}>0$, because the sign of the coefficient of $-dt^{2}$
\[
\beta (r,\theta)=1-\frac{2mr}{r^{2}+a^{2}\cos^{2}\theta}
\]
changes. Function $\beta (r,\theta)$ is null on the hypersurface
\[
r=m+\sqrt{m^{2}-a^{2}\cos^{2}\theta},
\]
and positive in the region $M^{a}$ outside this limit,
\begin{equation}\label{y}
M^{a}=\R\times\{x\in \R^{3}:r>m+\sqrt{m^{2}-a^{2}\cos^{2}\theta}\}.
\end{equation}
So, this region endowed with the metric (\ref{ep}) is stationary and is called {\it stationary Kerr spacetime}. Recall that if the rotating body covers the stationary limit hypersurface, then the gravitational field generated by the body is stationary (out of the body). 

Next, we  obtain an expression for the Hessian of a function as in Theorem \ref{t3}, applicable to study the convexity of stationary regions type
\begin{equation} \label{regmas}
M^{a}_{\epsilon}=\R\times\{x\in \R^{3}: r>m+\sqrt{m^{2}+\epsilon -a^{2}\cos^{2}\theta}\}\quad \epsilon >0,
\end{equation}
as in \cite{GM}, \cite{Ma}. Consider the function
\[
\phi_{a}(r,\theta)=\frac{1}{2}(r^{2}-2mr+a^{2}\cos^{2}\theta).
\]
Clearly

\[
\partial M^{a}_{\epsilon}=\R\times\{(r,\theta,\varphi): \phi_{a} (r,\theta)=\frac{1}{2}\epsilon \}.
\]
Since the radial component of the gradient of $\phi_{a}$ with respect to the Euclidean metric in $\R^{3}$ is equal to $r-m>0$, we have that $\partial M^{a}_{\epsilon}$ is smooth. 

We have not only that $\phi_{a}$ is independent of $t$ but also that is independent of $\varphi$ and, thus, $\delta (\phi_{a})=0$. So, as a consequence of Theorem \ref{t3}, % (\ref{e13})

\begin{corollary} For any  vector $\tilde{v}= t'_0\partial_t + v$ tangent to $M^a$,
\begin{equation}\label{eestr1}
{\rm Hess}_{\phi_{a}}(p)[\tilde{v},\tilde{v}]=  {\rm Hess}_{\phi_{a}}^{R}(x)[v,v]
+ t'_0{\rm rot}\delta(\nabla\phi_a,v) -\frac{1}{2}t'^{2}_{0}\nabla^{R}\beta(\phi_{a}). 
\end{equation}
\end{corollary}
Recall that $\partial M^{a}_{\epsilon}$ is (time, space or light) convex if 
${\rm Hess}_{\phi_{a}}(p)[\tilde{v},\tilde{v}]\leq 0$ for any (time, space or light) vector 
$\tilde{v}$ tangent to $\partial M^{a}_{\epsilon}$. This can be checked from (\ref{eestr1}) because the three terms in the right-hand side are directly computable. Indeed, 
let $\gamma(s)=(r(s),\theta (s),\varphi (s))$ be
a geodesic in $M_{0}$ such that $\gamma(0)=x$ and $\gamma'(0)=v (\equiv (r'_0,\theta'_0, \phi'_0)$). Putting $h(s)=\phi_{a}(\gamma(s))$ we obtain ${\rm Hess}_{\phi_{a}}^{R}(x)[v,v]=h''(0)$. On the other hand, since $\gamma(s)$ is a geodesic, %we have
%\[
%x_{1}''+\sum_{i,j=1}^{3}\Gamma_{i,j}^{1}x_{i}'x_{j}'=0 \quad x_{2}''+\sum_{i,j=1}^{3}\Gamma_{i,j}^{2}x_{i}'x_{j}'=0 
%\]
%($x_{1}\equiv r$, $x_{2}\equiv \theta$, $x_{3}\equiv \varphi$), which reduces to
\begin{equation}\label{ee}
\begin{array}{c}
r''=-\Gamma_{1,1}^{1}r'^{2}-2\Gamma_{1,2}^{1}r'\theta '-\Gamma_{2,2}^{1}\theta'^{2}-\Gamma_{3,3}^{1}\varphi '^{2} \\
  \theta''=-\Gamma_{1,1}^{2}r'^{2}-2\Gamma_{1,2}^{2}r'\theta '-\Gamma_{2,2}^{2}\theta'^{2}-\Gamma_{3,3}^{2}\varphi '^{2},
\end{array}
\end{equation}
where $\Gamma_{i,j}^{k}$ are the Christoffel simbols  for $<.,.>_{R}$.
Thus, replacing (\ref{ee}) in $h''(0)$ and taking into account that   
$r'_0=\frac{a^{2}\sin2\theta}{2(r-m)}\theta'_0$ (i.e., $\tilde{v}$ is tangent to $\partial M^{a}_{\epsilon}$):
\begin{equation}\label{eestr2}
\begin{array}{ll}
{\rm Hess}_{\phi_{a}}^{R}(x)[v,v]= & \theta_0 '^{2}\left(\frac{a^{4}\sin^{2}2\theta}{4(r-m)^{2}}-\Gamma_{1,1}^{1}\frac{a^{4}\sin^{2}2\theta}{4(r-m)}-\Gamma_{1,2}^{1}a^{2}\sin2\theta -\Gamma_{2,2}^{1}(r-m)\right. \\  & \left. -a^{2}\cos2\theta  +\Gamma_{1,1}^{2}\frac{a^{6}\sin^{3}2\theta}{8(r-m)^{2}}+\Gamma_{1,2}^{2}\frac{a^{4}\sin^{2}2\theta}{2(r-m)}+\Gamma_{2,2}^{2}\frac{a^{2}\sin 2\theta}{2}\right) \\ & +\varphi_0 '^{2}\left(-\Gamma_{3,3}^{1}(r-m)+\Gamma_{3,3}^{2}\frac{a^{2}\sin2\theta}{2}\right).
\end{array}
\end{equation}
%Now consider a geodesic $\gamma$ in $M^{a}_{\epsilon}$ with $\gamma(0)=p$, $\gamma'(0)=v$. Putting $q=<\gamma',\gamma'>$, $q_{0}=<\gamma'_{0},\gamma'_{0}>_{R}$ we have at $s=0$
%\[
%q_{0}=q-2g_{3,4}\varphi 't'-g_{4,4}t'^{2}.
%\]
%Thus, we obtain at $s=0$
%\[
%\theta '^{2}=\frac{4(r-m)^{2}}{\alpha(r,\theta)}q-\frac{4(r-m)^{2}g_{3,3}}{\alpha(r,\theta)}\varphi '^{2}-\frac{8(r-m)^{2}g_{3,4}}{\alpha(r,\theta)}\varphi 't'-\frac{4(r-m)^{2}g_{4,4}}{\alpha(r,\theta)}t'^{2}
%\]
%where $\alpha(r,\theta)=g_{1,1}a^{4}\sin^{2}2\theta+4g_{2,2}(r-m)^{2}$.
A straightforward computation shows:
\begin{equation} \label{eqeq}
\begin{array}{l}
\frac{1}{2}\nabla^{R}\beta (\phi_{a})=\overline{\Gamma}_{4,4}^{1}(r-m)-\overline{\Gamma}_{4,4}^{2}\frac{a^{2}\sin2\theta}{2} \\ 
{\rm rot}\delta(\nabla\phi_{a},v)=\left(\overline{\Gamma}_{3,4}^{2}a^{2}\sin2\theta - \overline{\Gamma}_{3,4}^{1}2(r-m)\right)\varphi'_{0}, 
\end{array}
\end{equation}
where $\overline{\Gamma}^{k}_{i,j}$ are the Christoffel simbols for $<.,.>$.
Summing up, substituting (\ref{eestr2}) and (\ref{eqeq}) in (\ref{eestr1}), a general expression for the Hessian of a vector  $\tilde{v} \equiv (t'_0, r'_0, \theta'_0, \phi'_0)$ tangent to 
$\partial M^{a}_{\epsilon}$ is obtained in terms of $t'_0,  \theta'_0, \phi'_0$.
From this expression,  one can  study when the region $M^{a}_{\epsilon}$ is (time, light or space) convex directly (compare with \cite[Ch. 7]{Ma}). 
%Finally, the subset $M^{a}$ is not spatially convex. This must be a direct consequence of the fact that fixed $p\in \partial M^{a}$
%\[(\ref{elarga})
%{\rm Hess}_{\phi_{a}}(p)[\tilde{v},\tilde{v}]={\rm Hess}_{\phi_{a}}^{R}(x)[\tilde{v},\tilde{v}]
%\]
%if $\tilde{v}\in T_{x}M_{0}\subseteq T_{p}M$.

\begin{corollary} \label{cec}
$M^{a}_{\epsilon}$ is not space convex for any $a^2\leq m^2, \epsilon >0$. 
\end{corollary}
{\it Proof.} 
From (\ref{eestr1}), if $\tilde{v}\in T_{p}M^{a}$ is tangent to $M_{0}$ then ${\rm Hess}_{\phi_{a}}(p)[\tilde{v},\tilde{v}]$ is the right-hand side of (\ref{eestr2}).
%\[
%{\rm Hess}_{\phi_{a}}(p)[\tilde{v},\tilde{v}]={\rm Hess}_{\phi_{a}}^{R}(x)[\tilde{v},\tilde{v}].
%\]
%But from  
%\begin{equation}
%\begin{array}{ll}
%{\rm Hess}_{\phi_{a}}^{R}(x)[\tilde{v},\tilde{v}]= &  \theta_0 '^{2}
%\left(\frac{a^{4}\sin^{2}2\theta}{4(r-m)^{2}}-\Gamma_{1,1}^{1}\frac{a^{4}\sin^{2}2\theta}{4(r-m)}-\Gamma_{1,2}^{1}a^{2}\sin2\theta -(r-m)\Gamma_{2,2}^{1} \right.\\
%& \left. -a^{2}cos2\theta  +\Gamma_{1,1}^{2}\frac{a^{6}\sin^{3}2\theta}{8(r-m)^{2}}+\Gamma_{1,2}^{2}\frac{a^{4}\sin^{2}2\theta}{2(r-m)}+\Gamma_{2,2}^{2}\frac{a^{2}\sin 2\theta}{2}\right) \\ & +\varphi_0 '^{2}\left(-\Gamma_{3,3}^{1}(r-m)+\Gamma_{3,3}^{2}\frac{a^{2}\sin2\theta}{2}\right)
%\end{array}
%\end{equation}
So, if $p=(t=0,r=m+\sqrt{m^{2}+\epsilon},\theta=\frac{\pi}{2},\varphi=0)\in \partial M^{a}_{\epsilon}$ and $\tilde{v}=\partial_{\theta}\in T_{p}\partial M^{a}_{\epsilon}$: 
%spacelike we obtain
\[
{\rm Hess}_{\phi_{a}}(p)[\tilde{v},\tilde{v}]=-(r-m)\Gamma^{1}_{2,2}+a^{2}=\frac{r(r-m)\triangle (r)}{\lambda(r,\theta)}+a^{2}>0. \quad \quad \Box
\]

\section{Non geodesic connectedness}\label{sad}

In this section we study the non geodesic connectedness of some regions of the slow ($a^{2}<m^{2}$), extreme ($a^{2}=m^{2}$) and fast ($a^{2}>m^{2}$) Kerr spacetime.
\begin{theorem}\label{p1}
Stationary Kerr spacetime $M^{a}$  with $0<a^2\leq m^2$ 
is not geodesically connected.
\end{theorem}
{\it Proof}. The first integrals of the geodesic equations of Kerr spacetime are 
\begin{equation}\label{eq}
\begin{array}{l}
\lambda(r,\theta)\varphi '=\frac{\D(\theta)}{\sin^{2}\theta}+a\frac{\P(r)}{\triangle (r)} \\ \lambda(r,\theta) t'=a\D(\theta)+(r^{2}+a^{2})\frac{\P(r)}{\triangle (r)} \\ \lambda(r,\theta)^{2}r'^{2}=\triangle (r) (qr^{2}-K)+\P^{2}(r) \\ \lambda(r,\theta)^{2}\theta '^{2}=K +qa^{2}\cos^{2}\theta-\frac{\D^{2}(\theta)}{\sin^{2}\theta}
\end{array}
\end{equation}
where
\[
\begin{array}{l}
\D(\theta)=L-Ea\sin^{2}\theta \\ \P(r)=(r^{2}+a^{2})E-La
\end{array}
\]
and being $q$ (normalization of the geodesic; rest mass), $K$ (Carter constant), $L$ (angular momentum) and $E$ (energy measured by observers in $\partial_t$) constants (we follow the notation in \cite[Chapter 4]{O}). If $\gamma(s)$ is a geodesic joining the points in the $z-$axis $p_{0}\equiv(t_{0}=0,r_{0},\theta_{0}=0)$ and $p_{1}\equiv(t_{1}=0,r_{1},\theta_{1}=\pi)$, $(r_{+}<)r_{0}<r_{1}$, in particular, $\gamma$ reaches the $z-$axis and, from the last equation in (\ref{eq}),  $L=0$ (otherwise $\lambda(r,\theta)^{2}\theta '^{2}$ would be negative near the z-axis). Then, the equation for $t$ reduces to
\[
\lambda(r,\theta) t'=E\left(\frac{(r^{2}+a^{2})^{2}}{(r-r_{-})(r-r_{+})}-a^{2}\sin^{2}\theta \right).
\]
% So, if we reparameterize this geodesic by $r$ we can write
%\[
%\triangle t=\int_{r_{0}}^{r_{1}}\frac{E[\frac{(r^{2}+a^{2})^{2}}{(r-r_{-})(r-r_{+})}-a^{2}\sin^{2}\theta ]}{\sqrt{(r-r_{-%})(r-r_{+})(qr^{2}-k)+[(r^{2}+a^{2})E]^{2}}}dr
%\]
%where the integral is considered in a general sense (possibly with rebounds).
But
\[
\frac{(r^{2}+a^{2})^{2}}{(r-r_{-})(r-r_{+})}-a^{2}\sin^{2}\theta
\]
is positive in $r\in (r_{+},+\infty)$ so if $\gamma$ satisfies $\triangle t=t_{1}-t_{0}=0$, necessarily $E=0$ and thus
\begin{equation}\label{ex}
\begin{array}{l}
\lambda(r,\theta)^{2}r'^{2}=\triangle (r) (qr^{2}-K) \\ \lambda(r,\theta)^{2}\theta '^{2}=K +qa^{2}\cos^{2}\theta.
\end{array}
\end{equation}

%\begin{equation}\label{ez}
%\triangle \theta =\int_{r_{0}}^{r_{1}}\frac{\sqrt{K+qa^{2}\cos^{2}\theta}}{\sqrt{(r-r_{-})(r-r_{+})(qr^{2}-K)}}dr.
%\end{equation}
Even more, from (\ref{ex})  we can ensure that $q>0$ and $0\leq\frac{K}{q}\leq r_{0}^{2}$ because, otherwise,
 either $\lambda(r,\theta)^{2}r'^{2}$ or $\lambda(r,\theta)^{2}\theta'^{2}$ would be negative at some point of $\gamma$. 
So, any zero $r^{*}$ of $\frac{\lambda(r,\theta)^{2}r'^{2}}{q}=(r-r_{-})(r-r_{+})(r^{2}-\frac{K}{q})$ cannot be greater than $r_{0}$. Now choose  $r_{1}<2m$. Then if, say $r(0)=r_{0}$, $r(1)=r_{1}$, at the point $s_{0}\in (0,1)$ such that $\theta (s_{0})=\frac{\pi}{2}$ we have $r(s_{0})>2m$ (see (\ref{y})) and $r'(s_{0})\neq 0$. So, if $r'(s_{0})>0$ (resp. $<0$) then $r(s_{0})<r(1)$ (resp. $r(0)>r(s_{0})$), in contradiction with $r_0, r_1 <2m$. $\Box$

Note that the previous proof can be extended in order to prove that (stationary or not) regions $R$ of outer Kerr spacetime with $a^{2}\leq m^{2}$ satisfying $r>r_{+}+\nu$ $(\nu >0)$ are not geodesically connected\footnote{Recall that this proof {\it cannot} be extended to the case $\nu=0$, which is geodesically connected \cite{FS-pr}.}. In fact, as $q>0$, $0\leq \frac{K}{q}\leq r_{0}^{2}$ and, so, any zero of $\frac{\lambda(r,\theta)^{2}r'^{2}}{q}$ is not greater than $r_{0}$, if we reparametrize $\gamma$ by $r$ then either 
\[
\triangle \theta =\int_{r_{0}}^{r_{1}}\frac{\sqrt{\frac{K}{q}+a^{2}\cos^{2}\theta}}{\sqrt{(r-r_{-})(r-r_{+})(r^{2}-\frac{K}{q})}}dr
\]
or, if $r'$ vanishes at a point $r^{*}$ ($<r_{0}$), perhaps:
\[
\triangle \theta =\int_{r^{*}}^{r_{0}}\frac{\sqrt{\frac{K}{q}+a^{2}\cos^{2}\theta}}{\sqrt{(r-r_{-})(r-r_{+})(r^{2}-\frac{K}{q})}}dr+\int_{r{*}}^{r_{1}}\frac{\sqrt{\frac{K}{q}+a^{2}\cos^{2}\theta}}{\sqrt{(r-r_{-})(r-r_{+})(r^{2}-\frac{K}{q})}}dr.
\]
As $\gamma$ must lie in $R$ then $r_{0},r_{1},r^{*}>r_{+}+\nu $; so, taking $r_{0},r_{1}$ close enough to $r_{+}+\nu$ we obtain necessarily $\triangle \theta$ small, which contradicts that $\theta_{1}-\theta_{0}=\pi$.

Moreover, no region $M^{a}_{\epsilon}$ ($a^{2}\leq m^{2}$), $\epsilon >0$ is  geodesically connected because of the following: 
(i) it lies in the region $r>m+\sqrt{m^{2}-a^{2}+\epsilon}$, and
(ii) the two points of this region non-connectable by geodesics found above, lie in
$M^{a}_{\epsilon}$. Summing up:

\begin{corollary}
Regions (stationary or not) of outer Kerr spacetime with $0\leq a^{2}\leq m^{2}$
 determined by $r>r_{+}+\nu$ for some $\nu >0$ are not 
geodesically connected.

Regions $M^{a}_{\epsilon}$ ($0\leq a^{2}\leq m^{2}$) are not geodesically connected 
for any $\epsilon >0$.
\end{corollary}
\begin{remark}\label{rrr} {\rm In the fast Kerr spacetime regions $M^{a}$ are again those with $\beta >0$ (compare with (\ref{y})) and regions $M^{a}_{\epsilon}$ do not have a natural sense. However, the same arguments work because $\triangle (r)=r^{2}-2mr+a^{2}$ admits a positive lower bound when $r>0$ which only depends on $a$. On the other hand, recall that in fast Kerr spacetime, one can consider $r\in\R$ and, thus to check the non-geodesic connectedness of $p_{0}=(t_{0},r_{0}<0,\theta_{0}=0)$ with $p_{1}=(t_{1},r_{1}>2m,\theta_{1}=\frac{\pi}{2},\varphi_{1})$, which lie in the stationary part. Summing up, we obtain}
\end{remark}
\begin{theorem}\label{te} (i) Stationary fast Kerr spacetime is not geodesically connected (if we assume $r>0$ as well as if $r\in \R$). 

(ii) Regions (stationary or not) of fast Kerr spacetime determined by $r>\nu$ for some $\nu>0$ are not geodesically connected.

(iii) The whole fast Kerr spacetime (including non-stationary regions and $r\in\R$) is not geodesically connected.
\end{theorem}

\section{Geodesic connectedness of Schwarzschild spacetime}

In this section we prove that given two points in $M^{a=0}$ there exist a geodesic joining them. Previously, we need the following technical result:
\begin{lemma}\label{l1} Let $\{f_{n}(x)\}_{n}$ be a sequence of continuous functions on $[a_{n},b]\subseteq \R$, $a_{n}\rightarrow a<b$ satisfying $0<c\leq f_{n}(x)\leq C$ for all $n$, and let $\{p_{n}(x)\}_{n}$ be  a sequence of polynomials with degree bounded in $n$ satisfying for all $n$: $p_{n}(a_{n})=0$, $p'_{n}(a_{n})=S_{n}>0$ and $p^{k)}_{n}(a_{n})\geq 0$ for $k\geq 2$.

(i) If $\{S_{n}\}_{n}\rightarrow \infty$, then
\[
\int_{a_{n}}^{b}\frac{f_{n}(x)}{\sqrt{p_{n}(x)}}dx\rightarrow 0.
\]

(ii) If $\{S_{n}\}_{n}\rightarrow 0$ and $p_{n}^{k)}(a_{n})$ admits an upper bound for  $k\geq 2$ and all $n$, then
\[
\int_{a_{n}}^{b}\frac{f_{n}(x)}{\sqrt{p_{n}(x)}}dx\rightarrow \infty.
\]
\end{lemma}

{\it Proof.} (i) Consider the sequence of polynomials $\{q_{n}(x)\}_{n}$, $q_{n}(x)=S_{n}(x-a_{n})\leq p_{n}(x)$ defined on $[a_{n},b]$. As $\int_{a}^{b}\frac{C}{\sqrt{x-a}}dx<\infty$ we have
\[
\int_{a_{n}}^{b}\frac{f_{n}(x)}{\sqrt{p_{n}(x)}}dx\leq \int_{a_{n}}^{b}\frac{C}{\sqrt{q_{n}(x)}}dx=\frac{1}{\sqrt{S_{n}}}\int_{a_{n}}^{b}\frac{C}{\sqrt{x-a_{n}}}dx \rightarrow 0.
\]

(ii) Because of the boundedness of $p^{k)}_{n}(a_{n})$, there exists $M>0$ such that $q_{n}(x)=S_{n}(x-a_{n})+M(x-a_{n})^{2}\geq p_{n}(x)$ on $[a_{n},b]$ for all $n$. Then
\begin{equation}\label{ccc}
\int_{a_{n}}^{b}\frac{f_{n}(x)}{\sqrt{p_{n}(x)}}dx\geq\int_{a_{n}}^{b}\frac{c}{\sqrt{q_{n}(x)}}dx=\int_{a_{n}}^{b}\frac{c}{\sqrt{S_{n}(x-a_{n})+M(x-a_{n})^{2}}}dx.
\end{equation}
But the sequence of last integrands converges uniformly on compact subsets of $(a,b]$ to the function $\frac{c}{\sqrt{M(x-a)^{2}}}$. Therefore, as $\int_{a}^{b}\frac{c}{\sqrt{M(x-a)^{2}}}=\infty$ we obtain that the limit in (\ref{ccc}) is $\infty$. $\Box$

\begin{theorem}\label{t0} Schwarzschild spacetime $M^{a=0}$ is geodesically connected.
\end{theorem}

{\it Proof.} Given two arbitrary points $p_{0}$ and $p_{1}$, the spherical symmetry of $M^{a=0}$ allows us to assume $p_{0}=(t_{0},r_{0},\varphi_{0},\theta_{0}=\frac{\pi}{2})$, $p_{1}=(t_{1},r_{1},\varphi_{1},\theta_{1}=\frac{\pi}{2})$, $t_{0}\leq t_{1}$; we can also assume $r_{0}\leq r_{1}$ (the modifications if $r_{0}>r_{1}$ are obvious). If we consider only geodesics $\gamma(s)$ on the ecuatorial plane $\theta\equiv \frac{\pi}{2}$, then its first integrals are obtained taking $a=0$, $K\equiv L^{2}$ and $\theta\equiv \frac{\pi}{2}$ in (\ref{eq}), that is: 
\begin{equation}\label{oo}
\begin{array}{l} r^{2}\varphi'=L \\ r^{2}t'=E\frac{r^{3}}{r-2m} \\ r^{4}r^{'2}=r(r-2m)(qr^{2}-L^{2})+r^{4}E^{2}.
\end{array}
\end{equation}
Notice also that if $t_{0}=t_{1}$ we can consider only geodesics with $E=0$; otherwise we can normalize $E=1$. Let $s(r)$ be the inverse function (where it exists) of $r(s)$ given by (\ref{oo}); using $r$ as parameter in the other two equations (\ref{oo}):
\begin{equation}\label{ooo}
\begin{array}{l}
\frac{d\varphi}{dr}=\epsilon \frac{L}{\sqrt{r(r-2m)(qr^{2}-L^{2})+r^{4}E^{2}}} \\ \frac{dt}{dr}=\epsilon\frac{Er^{3}}{(r-2m)\sqrt{r(r-2m)(qr^{2}-L^{2})+r^{4}E^{2}}}
\end{array}
\end{equation}
on a certain domain, being $\epsilon\in \{\pm 1\}$. If one consider geodesics with $r'\neq 0$ at any point, then the geodesic can be reparametrized by $r$ (recall $\mid s(r_{1})-s(r_{0})\mid <\infty$) and the increments $\triangle t$, $\triangle\varphi$ can be calculated integrating directly in (\ref{ooo}). Nevertheless, we are going to see that $p_{0}$, $p_{1}$ can be always joined with a geodesic such that $r'(s)$ vanishes exactly at one point $s^{*}$, and $r^{*}=r(s^{*})$ satisfies $2m<r^{*}<r_{0}$. Recall that the denominator in (\ref{ooo})
\begin{equation}\label{cb}
h(r)=r(r-2m)(qr^{2}-L^{2})+r^{4}E^{2}
\end{equation}
will vanish at $r^{*}$. As $r(s)$ will go from $r_{0}$ to $r^{*}$ then necessarily $h'(r^{*})>0$ (notice that this implies $\mid s(r^{*})-s(r_{0})\mid <\infty$). As later on $r(s)$ will go from $r^{*}$ to $r_{1}$, then $h(r)>0$ if $r^{*}<r<r_{1}$; we will consider geodesics with $h(r_{1})>0$ too. Summing up, it is sufficient to find constants $E$, $q$, $L^{2}$ as well as $r^{*}\in (2m,r_{0})$ such that the following relations (\ref{ggg}) and (\ref{pppp}) hold:
\begin{equation}\label{ggg}
\begin{array}{l} h(r^{*})=0,\quad h'(r^{*})>0,\quad\quad h(r)>0\quad\hbox{on}\quad (r^{*},r_{1}];
\end{array}
\end{equation}
putting
\begin{equation}\label{cbb}
\begin{array}{l}
\triangle t=\int^{r_{0}}_{r^{*}}\frac{Er^{3}}{(r-2m)\sqrt{h(r)}}dr + \int_{r^{*}}^{r_{1}}\frac{Er^{3}}{(r-2m)\sqrt{h(r)}}dr \\ \triangle \varphi =\int^{r_{0}}_{r^{*}}\frac{L}{\sqrt{h(r)}}dr + \int_{r^{*}}^{r_{1}}\frac{L}{\sqrt{h(r)}}dr,
\end{array}
\end{equation}
then
\begin{equation}\label{pppp}
\begin{array}{l}
\triangle t=t_{1}-t_{0} \\ \triangle \varphi=\varphi_{1}-\varphi_{0}+2k\pi
\end{array}
\end{equation}
for some integer $k$. Moreover, if $t_{0}=t_{1}$ we can fix $E=0$, if $t_{0}<t_{1}$ we fix $E=1$.

We will consider first the case $t_{0}<t_{1}$ and, thus, 
\begin{equation}\label{cbbb}
h(r)=r(r-2m)(qr^{2}-L^{2})+r^{4}.
\end{equation} 
If we look for $r^{*}$ such that $h(r^{*})=0$ and $h'(r^{*})=S >0$ then the following two relations for the constants $q$ and $L^{2}$ are obtained:

\begin{equation}\label{e'eqq}
\begin{array}{l}
qr^{*2}-L^{2}=-\frac{r^{*3}}{r^{*}-2m} \\ q=\frac{r^{*2}}{2(r^{*}-2m)^{2}}-\frac{3r^{*}}{2(r^{*}-2m)}+\frac{S}{2r^{*2}(r^{*}-2m)}.
\end{array}
\end{equation}
Taking into account the dependences of $q$ on $(r^{*}-2m)$ in (\ref{e'eqq}), there exist $r^{*}_{L}\in (2m,r_{0})$ near enough to $2m$ such that if $r^{*}\in (2m,r^{*}_{L}]$ then
\begin{equation}\label{e'eeq}
q>0\quad \forall S>0
\end{equation}
(note also that $L^{2}>0$ if $q>0$). Moreover, in order to apply Lemma \ref{l1},
\begin{equation}\label{l'}\begin{array}{l}
h^{2)}(r^{*})=6qr^{*}(r^{*}-2m)+2(qr^{*2}-L^{2})+(4q+12)r^{*2} \\ h^{3)}(r^{*})=6q(r^{*}-2m)+(18q+24)r^{*} \\ h^{4)}(r^{*})=24q+24.
\end{array}
\end{equation}
Clearly $h^{3)}(r^{*}),h^{4)}(r^{*})>0$ and, taking into account (\ref{e'eqq}) again, $h^{2)}(r^{*})>0$, so $h(r)>0$ if $r>r^{*}$. Summing up, it is sufficient to find an element of $A\equiv \{(r^{*},S):r^{*}\in(2m,r^{*}_{L}],S\in(0,\infty )\}$ such that the corresponding $(q,L^{2})$ given from (\ref{e'eqq})  and the function $h(r)\equiv h(r,q,L^{2})$ in (\ref{cbbb}) satify (\ref{pppp}) with $\triangle t$, $\triangle \varphi$ as in (\ref{cbb}) and $E=1$.

Fix $r^{*}\in (2m,r^{*}_{L}]$ and consider $\{(r^{*},S_{n})\}_{n}$, $\{S_{n}\}_{n}\rightarrow \infty$, then taking 
\[
f_{n}(r)\equiv\frac{r^{3}}{r-2m}\quad\hbox{and}\quad p_{n}(r)\equiv h_{n}(r)
\]
with $h_{n}\equiv h(r,q(r^{*},S_{n}),L^{2}(r^{*},S_{n}))$, hypotheses of Lemma \ref{l1} (i) clearly hold on the interval $[a_{n},b]=[r^{*},r_{1}]$. Therefore,
\[
(\triangle t)_{n}=\int^{r_{0}}_{r^{*}}\frac{f_{n}(r)}{\sqrt{p_{n}(r)}}dr + \int_{r^{*}}^{r_{1}}\frac{f_{n}(r)}{\sqrt{p_{n}(r)}}dr\rightarrow 0.
\]
This also holds if we take a sequence $\{r^{*}_{n}\}\rightarrow r^{*}$ and compute $(\triangle t)_{n}$ for $(r^{*}_{n},S_{n})$. Analogously, if we consider $\{S_{n}\}_{n}\rightarrow 0$ then, from (\ref{e'eqq}) and (\ref{l'}), $h^{k)}_{n}(r^{*})$ admits an upper bound for $k\geq 2$ and all $n$ thus, from Lemma \ref{l1} (ii), $(\triangle t)_{n}\rightarrow \infty$. In conclusion, given $\{\epsilon_{n}\}_{n}$, $\epsilon_{n}>0$, $\epsilon_{n}\searrow 0$ there exists $\{\delta_{n}\}_{n}$, $\delta_{n}>0$, $\delta_{n}\searrow 0$, such that
\begin{equation}\label{jjj}\begin{array}{l}
\triangle t(r^{*},S)<t_{1}-t_{0}\quad \quad \hbox{when} \quad \quad (r^{*},S)\in [2m+\epsilon_{n},r^{*}_{L}]\times [\frac{1}{\delta_{n}},\infty) \\ \triangle t(r^{*},S)>t_{1}-t_{0}\quad \quad \hbox{when} \quad \quad (r^{*},S)\in [2m+\epsilon_{n},r^{*}_{L}]\times (0,\delta_{n}]
\end{array}
\end{equation}
Next, we use topological arguments based on Brouwer's degree $deg$ (see the general viewpoint in \cite{L-S}). Essentially, we will prove that among the zeroes of $\triangle t-t_{1}+t_{0}$ given by (\ref{jjj}) there is a $(r^{*},S)\in A$ such that $\triangle \varphi$ satisfies (\ref{pppp}).
First, we prove
\begin{lemma}\label{mmm} There exists a connected subset ${\cal C}_{n}$ of zeroes of $\triangle t-t_{1}+t_{0}$ such that
\[
{\cal C}_{n}\cap (\{2m+\epsilon_{n}\}\times (\delta_{n},\frac{1}{\delta_{n}}))\neq \emptyset\quad \hbox{and} \quad {\cal C}_{n}\cap (\{r^{*}_{L}\}\times (\delta_{n},\frac{1}{\delta_{n}}))\neq \emptyset
\]
for every $n\in \N$. 
\end{lemma}
{\it Proof of Lemma \ref{mmm}.} Applying \cite[Lemma 3.4]{N} to the function
\[
\begin{array}{lcll}
 {\cal F}_{n}: & [2m+\epsilon_{n},r^{*}_{L}]\times (\delta_{n},\frac{1}{\delta_{n}}) & \rightarrow & X\equiv \R \\

           & (r^{*},S)       & \mapsto & \triangle t(r^{*},S)-t_{1}+t_{0}+S
\end{array}
\]
and using (\ref{jjj}), it is sufficient to prove (in the notation of \cite{N}):
\[
i_{X}({\cal F}_{n,r^{*}_{L}},G)\equiv deg(Id-{\cal F}_{n,r^{*}_{L}},G,0)\neq 0
\]
where ${\cal F}_{n,r^{*}_{L}}(S)={\cal F}_{n}(r^{*}_{L},S)$ and $G=(\delta_{n},\frac{1}{\delta_{n}})$ (recall that when $f\in C^{1}(a,b)\cap C[a,b]$, $f(a)\neq 0\neq f(b)$ and $f'(x)\neq 0$ if $f(x)=0$ then $deg(f,(a,b),0)=\Sigma_{x\in f^{-1}(0)}\hbox{sign}f'(x)$).
But the affine map
\[
\begin{array}{lcll}
 \hat{{\cal F}}_{n}: & (\delta_{n},\frac{1}{\delta_{n}}) & \rightarrow & \R \\

           & S       & \mapsto & 2\frac{1-\delta_{n}S}{1-\delta_{n}^{2}}-1+S
\end{array}
\]
has obviously $deg(Id-\hat{{\cal F}}_{n},G,0)=1$, and $deg(Id-{\cal F}_{n,r^{*}_{L}},G,0)=deg(Id-\hat{{\cal F}}_{n},G,0)$ (the map $\lambda \mapsto Id-{\cal F}_{n,r^{*}_{L}}+\lambda ({\cal F}_{n,r^{*}_{L}}-\hat{{\cal F}_{n}})$, $\lambda \in [0,1]$  is a homotopy from $Id-{\cal F}_{n,r^{*}_{L}}$ to $Id-\hat{{\cal F}}_{n}$   without zeroes on the boundary from (\ref{jjj})) which concludes the proof. $\Box$

Therefore, we obtain some $(r^{*}_{n},S_{n})\in {\cal C}_{n}$ with $r^{*}_{n}=2m+\epsilon_{n}$. Taking now in Lemma \ref{l1} (ii)
\[
f_{n}(r)\equiv 1\quad\hbox{and}\quad p_{n}(r)\equiv \frac{h_{n}(r)}{L_{n}^{2}},
\]
one checks from (\ref{e'eqq}) and (\ref{l'}) that its hypotheses hold on the intervals $[a_{n},b]=[r^{*}_{n},r_{1}]$, obtaining
\begin{equation}\label{uuu}
(\triangle \varphi)_{n}=\int^{r_{0}}_{r^{*}_{n}}\frac{f_{n}(r)}{\sqrt{p_{n}(r)}}dr + \int_{r^{*}_{n}}^{r_{1}}\frac{f_{n}(r)}{\sqrt{p_{n}(r)}}dr\rightarrow \infty.
\end{equation}
On the other hand, from (\ref{jjj}), the points in ${\cal C}_{n}$ with $r^{*}=r^{*}_{L}$ have $S\in (\delta_{1},\frac{1}{\delta_{1}})$, thus $L^{2}$ is upper bounded for these points and all $n$ (see (\ref{e'eqq})) and so is $\triangle \varphi$. This fact, (\ref{uuu}) and the connectedness of ${\cal C}_{n}$ imply the existence of $(r^{*},S)\in A$ such that (\ref{pppp}) holds, as required.

%Finally, consider the connected subset ${\cal C}=limsup_{n}\{ {\cal C}_{n} \}$ which satisfies
%\[
%{\cal C}\cap (\{r^{*}_{L}\}\times (0,\infty))\neq \emptyset
%\]
%and that given $n\in \N$ there exists $(r^{*}_{n},S_{n})\in {\cal C}$ such that $2m<r_{n}^{*}<2m+\frac{1}{n}$. Then, if $(r^{*}_{0},S_{0})\in {\cal C}\cap (\{r^{*}_{L}\}\times (0,\infty))$ again from Lemma \ref{l1} (i) we obtain that
%\[
%(\triangle \theta)_{n}-(\triangle \theta)_{0}\rightarrow \infty, \quad n\rightarrow \infty.
%\]  
%This joined with the connectedness of ${\cal C}$ implies the existence of $(r^{*},S) \in {\cal C}$ such that
%\[
%\begin{array}{l}
%\triangle t =t_{1}-t_{0} \\ 
%\triangle \varphi \in \{\varphi_{1}-\varphi_{0}+2(n-1)\pi : n\geq 1 \}.
%\end{array}
%\]

Finally, consider the case $t_{0}=t_{1}$ and thus,  put $E=0$. Now (\ref{cb}) becomes
\[
h(r)\equiv r(r-2m)(qr^{2}-L^{2}).
\]
By imposing $h(r^{*})=0$ and $h'(r^{*})=1$  we obtain the following values for the constants $q$ and $L^{2}$,
\begin{equation}\label{efqq}
\begin{array}{l}
qr^{*2}-L^{2}=0 \\ q=\frac{1}{2r^{*2}(r^{*}-2m)}.
\end{array}
\end{equation}
If we take $\{r^{*}_{n}\}_{n}\rightarrow 2m$, (\ref{efqq}) and formulas analogous to (\ref{l'}) imply that Lemma \ref{l1} (ii) can be applied to the functions 
\[
f_{n}(r)\equiv 1\quad\hbox{and}\quad p_{n}(r)\equiv \frac{h_{n}(r)}{L^{2}_{n}},
\]
on the intervals $[a_{n},b]=[r_{n}^{*},r_{1}]$. Thus, we obtain (\ref{uuu}) and, so, the existence of $r^{*}\in (2m,r_{0})$ such that (\ref{pppp}) holds. $\Box$

\begin{remark}\label{rr} {\rm The technique previously used in the proof of the geodesic connectedness of (outer) Schwarzschild spacetime is translatable to Schwarzschild black hole. In fact, now, we would use geodesics such that $r'(s)$ vanishes at $s^{*}$, with $r^{*}=r(s^{*})$ satisfying $(r_{0}\leq)r_{1}<r^{*}<2m$.}
\end{remark}

%So, assuming that $q>0$ we have
%\[
%\triangle \theta =\int_{r_{0}}^{r_{1}}\frac{\sqrt{\frac{K}{q}+a^{2}\cos^{2}\theta}}{\sqrt{(r-r_{-})(r-r_{+})(r^{2}-%\frac{K}{q})}}dr,\quad 0\leq \frac{K}{q}\leq r_{0}^{2}
%\]
%where the restriction on $\frac{K}{q}$ is imposed by requiring that both radicands are not negative. But if our %geodesic does not have  any rebound then $\triangle \theta \rightarrow 0<\theta_{1}-\theta_{0}=\pi$ when %$r_{0},r_{1}\rightarrow m+\sqrt{m^{2}+\epsilon (a)-a^{2}}$ independently of $\frac{K}{q}\in [0,r_{0}^{2}]$. %Finally, if my geodesic has a rebound $r^{*}$ to the left (a rebound to the right is not possible because %$h(r)=(r-r_{-})(r-r_{+})(r^{2}-\frac{K}{q})$ is strictly increasing) then $r^{*}\in (m+\sqrt{m^{2}+\epsilon (a)-%a^{2}},r_{0}]$ and, so, the derivatives of $h$ in $r^{*}$ different to zero have a positive inferior bound which %implies the same conclusion that in the previous case.

\end{document}